\documentstyle[sectioneqno,11pt]{article}
\newtheorem{theorem}{Theorem}[section]

  \title{ON A NEW ALGORITHM FOR THE  COMPUTATION OF ENCLOSURES FOR THE TITCHMARSH-WEYL $m$-FUNCTION \footnote{34E05, 34E20, 34L05, ODE's, Spectral theory,
Titchmarsh-Weyl $m$-function, verified computation}  
  }
\author{B.M. Brown*, M.S.P. Eastham, D.K.R. M$^c$Cormack, M. Plum}
\date{}
\begin{document}
\maketitle

\renewcommand{\baselinestretch}{1.5}
\large \normalsize     

 \newcommand{\beq}{\begin{equation}}
\newcommand{\enq}{\end{equation}}
\newcommand{\s}{ {\cal{x}}}
 \newcommand{\Ttildaa}{\overline{T} }
 \newcommand{\Ttilda}{\tilde{T}}

\section{Introduction}
Recently two of the present authors \cite{BrownEvansKirbyPlum1} reported on a
method for computing  safe bounds for the value of the Titchmarsh-Weyl $m$ function associated with the differential expression
\beq
My  \equiv \frac{1}{w}( -( py')' + q y) \label{eq:1.1}
\enq
defined over $[a,\infty), \; -\infty < a$, where $p,q,w$ are real-valued functions which
 satisfy $p^{-1},q,w\; \in L^1_{loc}[a,\infty)$  and $w(x) >0$ a.e. In the case that $w=1$ Weyl \cite{Weyl10} showed that  the differential equation
\beq
My=\lambda   y, \;\; \lambda \in   C_+ \cup C_-  \label{eq:1.2}
\enq
has at least one solution that belongs to the set
\begin{displaymath}
L^2_w(a,\infty) \equiv \{ f: \int_a^\infty  w\mid f \mid^2 dx < \infty \}.
 \end{displaymath}
The proof of this result introduced the Titchmarsh-Weyl $m$ function to the mathematics literature.
It was however  Titchmarsh   who investigated the properties of $m(\lambda)$ as an analytic function of $\lambda$ and established the connection between  the location of its poles, of necessity on the real line, and the eigenvalues of the differential equation (\ref{eq:1.2}).   
\par
The analytic form of the $m$ function is determined by    the  form of the $L^2$ solutions of
(\ref{eq:1.2}) and it is perhaps not suprising  that there are few examples of $m$ known in closed form.
For example if $a=0$ and $p=w=1$  then, when $q=\pm x^2$, the $m$ function is known as a rational function
of gamma functions, while if $q=\pm x$ it may be written in terms of Bessel functions. For  a detailed
discussion of the $m$ function, together with examples, see \cite{BE80}.
\par
In \cite{chaudhurieveritt} the central role  of the $m$ function in the spectral theory of (\ref{eq:1.1}) is established.
It is shown that the behaviour of the  $m$ function  near the real line classifies all   points of the real line   as belonging to the point spectrum, continuous spectrum, point-continuous spectrum or the resolvent set of  the self-adjoint operator  generated  in $L^2_w(a,\infty)$  by (\ref{eq:1.1}) together with initial conditions.  A further use for the $m$ function is in determining the best constant in Everitt's HELP inequality. See \cite{Everitt1} for further details.
\par
In view of the central  importance of the $m$ function and also in view of the difficulties  
in obtaining its values analytically, much effort has been expended in devising computational algorithms to 
estimate its value. These are reported on in a series of papers  which include \cite{BrownKirbyPryce1},\cite{BrownKirbyPryce2},\cite{BEEKa94}.
However these papers only contain estimates for the value of $m$ and do not seek to address the 
question  of absolute bounds on the error in the computations.
\par
 In \cite{BrownEvansKirbyPlum1} we reported on an algorithm to compute rigorous bounds for the $m$ function. This algorithm worked well 
for examples $q=x^\alpha,\;\; \alpha  \geq 2$  or $ \alpha =1$, but was shown to be computationally inefficient for 
examples  $q=-x^\alpha, \;   \alpha =1,2$, and the interval based algorithm needed in the computation  in \cite{BrownEvansKirbyPlum1} did not cover the case $q=\pm x^\alpha,\;\; 0 < \alpha < 2, \alpha \neq 1 $.  The purpose of this paper is to present two new algorithms which overcome these difficulties and enable the $m$ function now to be enclosed for a much wider class of problems.
\par
In section 2 we review the relevant extracts from the  theory of the Titchmarsh-Weyl $m$ function
that are needed to develop our algorithm.
Section 3 is devoted to recalling certain asymptotic results that are central to our method as well as
presenting an overview of  an interval based algorithm that is fundamental to the implementation of 
our $m$ computation.  Section 4 contains  the results of the numerical experiments that we 
have performed, while section 5 deals with the extension of the algorithm to overcome the problems encountered with  $q=\pm x ^\alpha, \; 0<\alpha <2,\; \alpha \neq 1$.
\section{Titchmarsh-Weyl limit-point, limit-circle classification}
In the  classical limit-point, limit-circle theory of (\ref{eq:1.2}) it is shown that,  starting from      a pair of solutions $\theta, \phi$ of (\ref{eq:1.2}) which  for strictly complex $\lambda$ satisfy
\begin{eqnarray}
\theta(a,\lambda)=0  & & ( p \theta')(a,\lambda)=1 \nonumber \\
\phi(a,\lambda)=-1  & & ( p \phi^{'})(a,\lambda)=0, \label{eq:2.1}
\end{eqnarray}
   there exists a complex-valued function $m(\lambda)$   such that 
\beq
\psi_0( \cdot, \lambda) = \theta(\cdot,\lambda) + m ( \lambda) \phi(\cdot,\lambda) \in L^2_w(a,\infty). \label{eq:2.2}
\enq
When, up to constant multiples,  there is precisely one solution  
of (\ref{eq:1.2}) in  $L^2_w( a,\infty)$,
we say (\ref{eq:1.1}) or (\ref{eq:1.2}) is limit-point at infinity. If all the solutions of (\ref{eq:1.2})
are in $L^2_w(a,\infty)$,  we say that  (\ref{eq:1.1}) or \mbox{ (\ref{eq:1.2})}  is limit-circle at infinity.  Further, the limit point, limit-circle classification  is determined  by $p,q,w$ and is independent of the strictly complex parameter $\lambda$. In this paper we shall be exclusively concerned with the limit-point case. The $m$  function is a Nevanlinna function, mapping the upper (lower) half-plane to itself, and as such has any singularities confined to the real line. From (\ref{eq:2.1}) and (\ref{eq:2.2})
it follows that
\beq
m(\lambda) = -\frac{\psi(a,\lambda)}{p(a) \psi'(a,\lambda)}  \label{eq:2.3}
\enq
where $\psi$ is any (non-zero) constant multiple of $\psi_0$.
The result (\ref{eq:2.3}) is the basis of our algorithm to compute $m(\lambda)$.
\par
We choose a point $X>0$ such that for $ x \in [X,\infty)$ we may develop an asymptotic expansion 
for $\psi(x,\lambda)$, together with a precise estimate on the error  committed.   This   expansion  enables us to determine   intervals in which $\psi(X,\lambda)$ and   $ \psi'(X,\lambda)$ lie, thus providing    
 initial data to an interval based  initial value solver that is used to compute    complex intervals which enclose
$\psi(a,\lambda)$ and  $ \psi'(a,\lambda)$.  The result (\ref{eq:2.3}) yields an interval which encloses $m(\lambda)$. In section 3 we review the asymptotic method and interval  ODE solver that is used to
perform these tasks.
\section{Overview of the components of the algorithm}
\subsection{Asymptotic theory}
The method that we   use to obtain   the   asymptotic  solution   (\ref{eq:1.2})
as $x \rightarrow \infty$ is the repeated diagonalization method of Eastham which is fully  explained in the book \cite{MSPE89} and here we shall be brief.
The method is concerned with  estimating and improving  error terms  in the asymptotic solution of 
the linear differential system
\beq
Z'(x)  = \rho(x) \{ D + R(x) \} Z(x) \;\;\;(a \leq x < \infty) \label{eq:3.1}
\enq
where $Z$ is an $n-$component vector, $\rho$ is a real or complex scalar factor, $D$ is a constant diagonal matrix
\begin{displaymath}
D={\rm dg } ( d_1,d_2,...,d_n)
\end{displaymath}
with distinct $d_k$ and $R$ is a perturbation such that 
\beq
R(x)=O(x^{-\delta}) \;\;\; ( x \rightarrow \infty) \label{eq:3.2}
\enq
for some $\delta >0$.
\par
If it is the case that $\rho R \in L(a,\infty)$, the Levinson asymptotic theorem can be applied to 
(\ref{eq:3.1}) to give solutions
\beq
Z_k(x) = \{ e_k + \eta_k(x) \} \exp ( d_k \int_a^x \rho (t) dt ) \label{eq:3.3}
\enq
where $e_k$ is the unit coordinate vector in the $k-$direction and  $\eta_k(x) \rightarrow 0$ as $ x \rightarrow \infty$. The size of the error term is related to the size  of $R$ as $x \rightarrow  \infty$, and therefore the accuracy of (\ref{eq:3.3}) can be improved if the perturbation $R$ can be reduced as $x \rightarrow \infty$. Under suitable conditions on $\rho$ and $R$ this improvement can be achieved by  a sequence of repeated transformations which lead to a  computational procedure to estimate the 
solution of  (\ref{eq:1.2}) together with a bound on  the associated  error.
\par
 The sequence of transformations may be obtained either by an {\it exact diagonalization} or
by an {\it approximate diagonalization} procedure.  These methods are discussed in detail in \cite{BEEM95}.
The exact diagonalization method  involves the explicit construction of an $n \times n $ matrix $T$ such that
\begin{displaymath}
T^{-1}(x) \{ D + R(x) \} T(x) = D_1(x)
\end{displaymath}
and this in turn requires the explicit eigenvectors of  $D+R(x)$ which, although available 
for the second order system, are not generally known for the $n-$th order system.
In this  investigation we choose to work with the more generally applicable  {\it approximate diagonalization 
method} which may be used for $n-$th order systems of differential equations.  A discussion on the asymptotic method of    exact diagonalization  as applied to estimating the $m$ function can be found in \cite{BEEM95}.
\subsubsection{Approximate diagonalization}
We   assume that $R$ is a differentiable  $n \times n $ matrix satisfying (\ref{eq:3.2}),  and  we
  define an $n \times n $ matrix $P$ by
\begin{displaymath}
PD-DP=R-{\rm dg}R
\end{displaymath}
with diagonal  entries  $p_{ii}=0$ and other entries
\beq
p_{ij}=r_{ij}/(d_j-d_i).  \label{eq:3.3a}
\enq 
We note that
$P=O(R)=O(x^{-\delta})$ and the construction of the $P$ matrix cancels out the dominant terms in the 
following  system (\ref{eq:3.5}). With 
$Z=(I+P)W
$,  (\ref{eq:3.1}) is transformed into  the system
\par
\beq
W'=\rho ( \tilde{D}+S)W
\label{eq:3.5}
\enq
where we have written  the $n \times n$ matrices
\begin{eqnarray}
\tilde{D}&=&D + {\rm dg } R \nonumber \\
 S=(I+P)^{-1}(RP -P {\rm dg }R-\rho^{-1}P')&=&O(x^{-2\delta})+O(\rho^{-1}x^{-\delta-1})
\label{eq:3.6}
\end{eqnarray}
which is of smaller magnitude than $R$. We write
\beq
(I+P)^{-1}=I-P+P^2+...+(-1)^ \nu P^\nu +(-1)^{\nu+1}(I+P)^{-1}P^{\nu+1} \label{eq:3.7}
\enq
and specify an order of magnitude
$O(x^{-K})
$ which we wish to achieve as an error in the asymptotic solution of  (\ref{eq:3.1}). Substituting (\ref{eq:3.7}) into (\ref{eq:3.6}), we have
\beq
S=V_2+...V_{M-1} + E \label{eq:3.8} 
\enq
where
\begin{displaymath}
V_m=O(x^{-m \delta})
\end{displaymath}
\begin{displaymath}
E=O(x^{-K})
\end{displaymath}
and
$(M-1) \delta < K \leq M \delta$.
Thus $\nu$ is chosen so that $P^{\nu+1}$ in (\ref{eq:3.7}) gives rise to terms which contribute to $E$ by at most this order of magnitude, and will be estimated in the final stages of the algorithm.
\par
This transformation procedure may be repeated for the $W$ system (\ref{eq:3.5}), but with a new $P$ defined in terms of $V_2$ which replaces $R$ in (\ref{eq:3.1}). We continue to use the matrix $D$ and not $\tilde{D}$
to simplify  the construction of an efficient computational algorithm. However this procedure  introduces additional terms into the analysis which must be eliminated at subsequent iterations of the algorithm.  A repetition of the above process leads to a new matrix $S$ {\it viz.}
\begin{displaymath}
S=V_3+ ...+V_{M-1} + E
\end{displaymath}
with new $V$'s and a new $E$. 
\par
The above ideas may be used to form the basis of an iterative procedure,  which can be implemented in the symbolic algebra system {\it Mathematica}, to compute the asymptotic
solutions of  (\ref{eq:3.1}).  Taking \mbox{(\ref{eq:3.1})} as a starting point with $m=1$, we have 
at the $m-$th stage 
\begin{eqnarray*}
Z'_m&=&\rho(D_m+R_m)Z_m, \\
R_m&=&V_{1m}+V_{2m}+...+V_{M-m,m}+E_m, \\
V_{jm}&=&O(x^{-(m+j-1)\delta}), \\
{\rm dg }V_{1m} &=& 0, \\
  D_m&=&D+\Delta_m, \\
Z_m&=&(I+P_m)Z_{m+1} 
\end{eqnarray*}
where $P_m$ is defined explicitly in terms of  $V_{1m}$ and $D$ as in (\ref{eq:3.3a}).
Thus $V_{1m}$ is eliminated at this stage.
\par
At the end of the process all the $V$'s are eliminated and this gives
\begin{displaymath}
Z_M' = \rho(D_M+R_M)Z_M 
\end{displaymath}
where
\begin{displaymath}
R_M = E_M = O(x^{-K})
\end{displaymath}
and  $M \delta \geq K $.
 The Levinson theorem   then yields the solution of the $Z_M$ system, and  reversal of the $M-1$ transformations gives
 \begin{displaymath}
Z= \{ \Pi_1^{M-1}(I+P_m)\} ( e + \eta ) \exp(\int_a^x \rho(...)dt)
\end{displaymath} 
with $\eta=O(x^{-K})$, ( see \cite{BEEM95} for further details).
\subsection{Interval ODE solver}
In this sub-section we introduce briefly  the concepts of interval arithmetic that we   need
to give a short account of Lohner's AWA algorithm.  For  an in-depth discussion of interval arithmetic, see \cite{AH83}, while Lohner's AWA algorithm is discussed in \cite{Lohphd} and \cite{BrownEvansKirbyPlum1}.
\par
Denoting any of the four basic arithmetic operations by  $\star$, we define,
for real intervals $[a],[b]$,  
\begin{displaymath}
[a] \star [b] = \{ a \star b \mid a \in [a],\; b \in [b] \}.
\end{displaymath}
Thus we can compute  an enclosure for $[a ]\star[ b]$ by  obtaining   computable upper, and lower bound, for $[a] \star[ b]$ which is derived from the lower and upper bounds of $[a],[b]$ respectively, by some directed rounding facilities.
Any algorithm that is realised on a computer consists of finitely many operations $\star$ and thus 
an enclosure for the  results of arithmetic operations  which constitute the algorithm may be computed. In practice this simple approach would soon lead to an explosion of the interval width but many sophisticated techniques are available to control this phenomenon \cite{AH83}.
\par
Lohner's approach to computing an enclosure of the solution of initial value problems
 is based on  the well known Taylor method for solving initial value problems.
Suppose  that a solution of the IVP
 \beq
u^{'}=f(x,u), \;\;u(0)=u_0, \label{eq:3.8a}
\enq
where $f:[0,\infty)  \times R^n \rightarrow R^n$ is sufficiently smooth, 
is known at some point $x_0$.  Then the solution at $x_0+h$ is
\beq
u(x_0+h)= u(x_0) + h \phi( x_0,h) + z_{x_0+h} \label{eq:3.9}
\enq
where  $  u(x_0) + h \phi( x_0,h) $ is the $(r-1)$--th degree Taylor  polynomial  of $u$ expanded about $x_0$ and  $z_{x_0+h}$ is the associated local error.
This method lends itself well to computation since the coefficients of the polynomial may be computed  via an automatic differentiation package by differentiating the 
differential equation (\ref{eq:3.8a}).  However the error term is not known exactly since   the standard   formulae give, for some unknown $\tau$,
\beq
z_{x_0+h}=u^{(r)}(\tau)h^{r}/r!, \;\;\; \tau \in [x_0,x_0+h].\label{eq:3.10}
\enq
In Lohner's algorithm, (\ref{eq:3.9}) is used to advance an enclosure  $[u(x_0)]$  for the solution
$u$ at $x_0$, to one for the solution $u$ at $x_0+h$ which we denote by   $[u(x_0+h)]$. A suitable  enclosure  for the error (\ref{eq:3.10})
is
\begin{displaymath}
[z_{x_0+h}]=f^{(r)} ( [ x_0,x_0+h],[u])h^r/{r!}
\end{displaymath}
provided  that an enclosure $[u]$ for $\{ u(x) : x_0 \leq x \leq x_0+h \}$
can be computed. This is achieved by the following means.
Choose some interval $[u^0] \supset [u(x_0)]$ and try to prove that 
\begin{displaymath}
[u]=[u(x_0)] + [0,h] \cdot f([ x_0,x_0+h],[u^0]) \subset [u^0].
\end{displaymath}
If this is true then Banach's fixed-point theorem implies  that $[u]$ is an enclosure for $u(x )$ for all $x \in ( x_0,x_0+h)$.
In order to achieve efficient  performance and tight bounds,  the   details 
of the algorithm  are more complex than  this short overview can show. We refer the reader to \cite{Lohphd} and \cite{BrownEvansKirbyPlum1} for a complete discussion of the method.
 \section{Results for $q=-x^\alpha,\;\;\alpha  =\;1,\; 2$}
In this section we discuss the computation of $m$ when $a=0$ and  $p=w=1$ and   the potential  $q= -x^\alpha,\;0 \leq  \alpha \leq 2 $. In terms of the asymptotic analysis presented in section 3 this means that we take $\delta = \alpha/2$. 
However, while the general algorithm is applicable to all $\alpha$ in this range, the implementation of 
Lohner's AWA interval ODE solver requires at least two derivatives  of the function $q$,
 see (\ref{eq:1.1}),
to be available at $x=0$. Clearly this is not possible for $0  < \alpha <2,\; \alpha \neq 1 $, and for $\alpha $  in this range a revised algorithm    is presented   in section 5.  Here we present an algorithm  to compute $m$ when $q=-x$ or $ q=-x^2$. We remark that the algorithm which we present here
is  also applicable to problems where $q=x^\alpha,  \; 2 \leq \alpha$ or $ \alpha =1$, while that in section 5 covers the case $0 < \alpha < 2$. 
\par
We first write (\ref{eq:1.2})  as the system 
\begin{displaymath}
\left ( \begin{array}{c} y_1 \\ y_2 \end{array} \right )' = 
\left ( \begin{array}{cc} 0 & 1  \\ q-\lambda & 0  \end{array} \right )
\left ( \begin{array}{c} y_1 \\ y_2 \end{array} \right )
\end{displaymath}
and as in \cite[chapter 2]{MSPE89} introduce the transformation 
\begin{displaymath}
T=   \left ( \begin{array}{cc}
1 & 1 \\
\sqrt{ q-\lambda} & -\sqrt{ q-\lambda}  \\
\end{array}
\right ) .
\end{displaymath} 
This  enables us to write (\ref{eq:1.2}) in the form (\ref{eq:3.1})  with   
\begin{displaymath}
\rho= \sqrt{q-\lambda},\;\;D= {\rm dg }(1,-1)
\end{displaymath}
and 
\begin{displaymath}
R=
\frac{q^{'}}{4(q-\lambda)^{3/2}} \left (
\begin{array}{cc}
-1 & 1\\
1 &-1 \\
\end{array}
\right ).
\end{displaymath}
We next apply $6$ iterations of the asymptotic algorithm of section 3.1 to obtain  bounds on   $\psi(x,\lambda)$ and
$\psi'(x,\lambda)$  $ (X \leq x <\infty)$ the $L^2[a,\infty)$ solution obtained from the asymptotic algorithm.  This gives intervals which enclose $\psi(X,\lambda)$
and $\psi'(X,\lambda)$ which  are   the initial data required by the    AWA algorithm. 
The asymptotic analysis and estimation of the error is performed  using   purpose written Mathematica code, the  detail of which is fully reported on in \cite{BEEM95}.  
\par	
The C-XSC implementation of Lohner's algorithm is used with purpose written shell script to interface the asymptotic results to the interval arithmetic code.

\begin{table}
\small
\begin{center}
\begin{tabular}{||c|c||} \hline
$\lambda$ & $m(\lambda)$ \\ \hline
\hline
-1 + i 	&   $ 0.723783^{7644403761}_{6822580120} + 0.428707^{1190558884}_{0519233491}$ i  \\ \hline
i 	&  $ 0.555050^{3090709121}_{2364882611} + 0.665360^{7324705136}_{6535858522} $ i  \\ \hline
0.5 + i &  $ 0.4240432^{898149527}_{260330653} + 0.700941^{1148568185}_{0374761820} $ i \\ \hline
0.1 + 0.1 i &  $ 0.59997^{30075913471}_{29027001456} + 1.09669^{40756787737}_{39419214239} $ i \\ \hline
1 + i 	& $0.3127262^{773973635}_{385936455}$ + $0.6886661^{844509032}_{343957574}$ i   \\ \hline
10 + 10 i & $ 0.101753^{7139152363}_{6973947839} + 0.2444418^{474486887}_{236960981} $ i  \\ \hline
1 + 0.5 i & $ 0.265972^{3056163538}_{2659812543} + 0.7947538^{875505129}_{307529690} $ i  \\ \hline
1 + 0.1 i & $ 0.187356^{8032104862}_{7636023492} + 0.885754^{9186094842}_{8512935614} $ i \\ \hline
 1 + 0.01 i & $ 0.1628908^{559868685}_{161172176} + 0.904793^{5039573860}_{4316287238} $ i  \\ \hline
 1 + 0.001 i &  $ 0.1602956^{820262119}_{421095763} + 0.906631^{1561255695}_{0831753883} $ i \\ \hline
 1 + 0.0001 i & $ 0.1600346^{705614590}_{306398719} + 0.906814^{1685974984}_{0955836808} $ i  \\ \hline
 1 + 0.00001 i & $ 0.1600085^{544702433}_{145481948} + 0.906832^{4622237933}_{3892031914} $ i  \\ \hline
\hline

\end{tabular}
\end{center}
\label{tab:4.1}
\caption{$X = 10$ and $\alpha = 1$ , where $\epsilon_6(10) = 1.09576673\!\times\!10^{-8}$. }
\end{table}
 \begin{table}
\small
\begin{center}
\begin{tabular}{||c|c||} \hline
$\lambda$ & $m(\lambda)$ \\ \hline
\hline
-1 + i 		& $ 0.7215463^{224474854}_{158603078} + 0.3676480^{842327160}_{790723077} $ i  \\ \hline
i 		& $ 0.6266570^{722135664}_{651135956} + 0.6266570^{722008661}_{651302206} $ i  \\ \hline
0.5 + i 	& $ 0.4999306^{961394431}_{890513378} + 0.7052934^{170794993}_{088942980} $ i \\ \hline
0.1 + 0.1 i 	& $ 0.8975088^{969388575}_{867879163} + 1.0328959^{633792185}_{527567132} $ i \\ \hline
1 + i 		& $ 0.367648^{0840684596}_{0792256377} +  0.7215463^{223138735}_{159732310} $ i   \\ \hline
10 + 10 i 	& $ 0.10206645^{54130719}_{40537547} + 0.24554208^{42914867}_{23862309} $ i  \\ \hline
1 + 0.5 i 	& $ 0.3372315^{718010392}_{665134664} + 0.86795375^{85127644}_{05598003} $ i  \\ \hline
1 + 0.1 i 	& $ 0.25612^{60389488290}_{59314609572} + 1.022821^{6685241975}_{5104881830} $ i \\ \hline
 1 + 0.01 i 	& $ 0.22527283^{50846718}_{03204735} + 1.061093^{7010720598}_{6911651063} $ i  \\ \hline
 1 + 0.001 i 	& $ 0.22185279^{81363766}_{33871469} + 1.0649430^{570902692}_{471166850} $ i \\ \hline
 1 + 0.0001 i 	& $ 0.2215072^{705977203}_{658584275} + 1.0653280^{561152201}_{461491257} $ i  \\ \hline
 1 + 0.00001 i 	& $ 0.2214726^{824283285}_{776892570} + 1.0653665^{564746155}_{465145155} $ i  \\ \hline
\hline

\end{tabular}
\end{center}
\label{tab:4.2}
\caption{$X = 10$ and $\alpha = 2$ , where $\epsilon_6(10) = 1.20325893\!\times\!10^{-9}$. }
\end{table}
 
\begin{table}
\small
\begin{center}
\begin{tabular}{||c|c||} \hline
$\lambda$ & $m(\lambda)$ \\ \hline
\hline
-1 + i 		& $ 0.723783723^{4729994}_{2172666} +  0.428707085^{5961411}_{3870145}$ i  \\ \hline
i 		& $ 0.5550^{505827338581}_{499628151809} + 0.66536^{09723872275}_{04136894624} $ i  \\ \hline
0.5 + i 	& $ 0.42404325^{87364205}_{71105477} + 0.70094107^{71574380}_{51492617} $ i \\ \hline
0.1 + 0.1 i 	& $ 0.^{6000855062924143}_{5998604127881420} + 1.096^{7992302211012}_{5887988999492} $ i \\ \hline
1 + i 	  	& $ 0.31272625^{81969598}_{78116266} + 0.688666159^{8412076}_{0256750} $ i   \\ \hline
10 + 10 i 	& $ 0.1017537056^{555976}_{544222} + 0.24444183557^{31793}_{16071} $ i  \\ \hline
1 + 0.5 i 	& $ 0.26597^{34567852943}_{11149849902} + 0.79475^{53834479725}_{23351057398} $ i  \\ \hline
1 + 0.1 i 	& $ 0.187356^{8022703960}_{7653048066} + 0.885754^{9051312494}_{8667747971} $ i \\ \hline
1 + 0.01 i 	& $ 0.162890^{9120955266}_{7609618708} + 0.904793^{5903178042}_{3484631105} $ i  \\ \hline

1 + 0.001 i 	& $ 0.1602956^{663626111}_{587434146} + 0.906631^{1260873454}_{1165602242} $ i \\ \hline

1 + 0.0001 i 	& $ 0.160034^{8624420599}_{4397312929} + 0.90681^{45096602290}_{37578832769} $ i  \\ \hline
1 + 0.00001 i 	& $ 0.160008^{5817370011}_{4882534979} + 0.906832^{5038096614}_{3509810573} $ i  \\ \hline
\hline

\end{tabular}
\end{center}
\label{tab:4.3}
\caption{$X = 40$ and $\alpha = 1$ , where $\epsilon_6(40) = 5.21570339\!\times\!10^{-15}$. }
\end{table}

In all cases the enclosures that we obtain for $m(\lambda)$ are in agreement with the closed form  of $m(\cdot)$ given in terms of either gamma functions or Bessel functions \cite{ECT62}, and evaluated by  numerical routines, (see   \cite{TB87} and \cite{Kuk72} for further
details). We further remark that the algorithm reported on in \cite{BEEK94}
could not perform the computation required to produce the above results. 

\newpage

 \section{Results for $q=-x^\alpha$,$0  < \alpha < 2$, $a=0,p=w=1$}
As we mentioned at the beginning of section 4 
 we have to  modify our algorithm to
deal with values of  $\alpha$  other than $1$ and $2$.  
We do this  by choosing some number $\epsilon >0$.  Then, instead of  solving  
 (\ref{eq:1.2})  over  the whole of $[0,X]$, we now  solve the equation  over the interval $[\epsilon,X]$. The following   Theorem 5.1 enables us to    obtain
an enclosure for $y(0)$ and  $y^{'}(0)$ in terms of the enclosures for  $y(\epsilon)$ and $ y^{'}(\epsilon)$.
\begin{theorem}
Let $c(x)=q(x)-\lambda$ and, for some $\epsilon >0$,  let $f,g \in C_1[0,\epsilon]$ satisfy $f(\epsilon)=y(\epsilon)$, $g(\epsilon)=y^{'}(\epsilon)$.   In addition let
\begin{displaymath}
b=\epsilon \int_0^\epsilon \mid c \mid dt < 1.
\end{displaymath}

Then
\begin{enumerate}
\item
\beq
\mid y(0)-f(0) \mid \leq \frac{1}{1-b} [ \int_0^\epsilon \mid f^{'} -g \mid dt + \epsilon \int_0^\epsilon
\mid g^{'}-cf \mid dt ]
\label{eq:5.1}
\enq
\item
 \beq
\mid y^{'}(0)-g(0) \mid \leq \frac{1}{1-b} [  \int_0^\epsilon\mid c \mid dt  \int_0^\epsilon \mid f^{'} -g \mid dt +  \int_0^\epsilon
\mid g^{'}-cf \mid dt ].
\label{eq:5.2}
\enq
\end{enumerate}
\end{theorem}
{\bf Proof:}
Define
\begin{displaymath}
u(x)= \left (\begin{array}{c} y(x) \\ y^{'}(x) \\\end{array} \right )
-\left( \begin{array}{c} f(x) \\ g(x)\\ \end{array} \right ).
\end{displaymath}
Then $u(\epsilon)=0$ and 
\begin{eqnarray*}
u^{'}(x)  &=&  \left (\begin{array}{c} y^{'}(x) \\c(x)y(x) \\\end{array} \right )
-\left (\begin{array}{c} f^{'}(x) \\ g^{'}(x) \\\end{array} \right  )\\
&=& \left ( \begin{array}{cc}0&1\\ c(x)&0\\\end{array}\right ) u(x) -\tau(x),
\end{eqnarray*}
where  
\begin{displaymath}
\tau= 
\left ( \begin{array}{c}f^{'}- g\\ g^{'}-cf\\\end{array}\right ).
\end{displaymath}
\par
 Define 
\begin{displaymath}
T:(C[0,\epsilon])^2 \rightarrow  (C[0,\epsilon])^2
\end{displaymath}
by
\begin{displaymath}
(Tv)(x) := \int_\epsilon^x 
\left (
\begin{array}{cc} 0&1 \\ c(t) & 0 \end{array} \right ) v(t) dt - \int_\epsilon^x \tau(t)dt,
\end{displaymath}
for $x \in [0,\epsilon ]$ and   $ v \in (C[0,\epsilon])^2$.
Then it follows that
\begin{displaymath}
u=Tu.
\end{displaymath}
We shall prove that $T$ has a globally unique fixed point in
\begin{displaymath}
U:= \{ v \in (C[0,\epsilon])^2: \mid v_1(x) \mid \leq \alpha_1, \; \mid v_2(x)\mid \leq \alpha_2 \;
(0 \leq x \leq \epsilon) \}
\end{displaymath}
where $\alpha_1$ and  $ \alpha_2$ denote the right-hand sides of the inequalities (\ref{eq:5.1}) and (\ref{eq:5.2}), respectively. To show this  we use the Banach  fixed point theorem.
This requires us to show
\begin{enumerate}
\item
$T$ is a contraction with respect to some suitable norm;
\item
$TU \subset U$.
\end{enumerate}
In order to show (1) let 
\begin{displaymath}
\parallel v \parallel := {\max}_{x \in [0,\epsilon]}\{ \max \{ \frac{1}{\sqrt{\epsilon}}
\mid v_1(x) \mid , \frac{1}{\sqrt{ \int_0^\epsilon \mid c \mid dt }} \mid v_2(x) \mid \}\}.
\end{displaymath}
Since for all $v,\tilde{v} \in (C[0,\epsilon])^2$ and $x \in [0,\epsilon]$,
\begin{eqnarray*}
\mid( Tv-T\tilde{v})_1(x) \mid  &\;\leq \;&  \int_x^\epsilon \mid v_2(t)-\tilde{v}_2(t) \mid dt  \leq  \epsilon
\parallel v-\tilde{v} \parallel \sqrt { \int_0^\epsilon \mid c \mid dt },  \\
\mid( Tv-T\tilde{v})_2(x) \mid  &\;\leq\; & \int_x^\epsilon \mid c(t) \mid \mid v_1(t)-\tilde{v}_1(t) \mid dt  \leq  \sqrt{\epsilon}
\parallel v-\tilde{v} \parallel  \int_0^\epsilon \mid c \mid dt . 
\end{eqnarray*}
we have
\begin{displaymath}
\parallel Tv -T \tilde{v} \parallel \leq \sqrt{ \epsilon \int_0^\epsilon \mid c \mid dt } \parallel v-\tilde{v}
\parallel.
\end{displaymath}
Since 
\begin{displaymath}
\epsilon \int_0^\epsilon \mid c \mid dt <1
\end{displaymath}
it follows that $T$ is a contraction.
\par
In order to establish part 2 we take $v \in U$. Then 
\begin{eqnarray*}
\mid (Tv)_1(x) \mid &\leq& \int_x^\epsilon \mid v_2 \mid dt + \int_x^\epsilon \mid \tau_1 \mid dt \\
&\leq& \epsilon \alpha_2 + \int_0^\epsilon \mid \tau_1 \mid dt \\
&=& \alpha_1
\end{eqnarray*}
and
\begin{eqnarray*}
\mid (Tv)_2(x) \mid &\leq& \int_x^\epsilon\mid c \mid  \mid v_1\mid dt + \int_x^\epsilon \mid \tau_2 \mid dt \\
&\leq&  \alpha_1 \int_0^\epsilon \mid c \mid dt  + \int_0^\epsilon \mid \tau_2 \mid dt \\
&=& \alpha_2.
\end{eqnarray*}
 Thus   $T$ has a  globally unique fixed point in $U$ which implies $u \in U$. In particular, 
\begin{displaymath}
\mid u_1(0) \mid \leq \alpha_1, \;\;\mid u_2(0) \mid \leq \alpha_2
\end{displaymath}
which establishes the theorem.
 \subsection{Numerical examples}
We   now turn  to two examples which illustrate the use of  Theorem 5.1
to compute enclosures for $m$.
These examples are: $p = w = 1$,  and
\begin{enumerate}
\item
$q(x)=-\sqrt{x}$;
\item
$q(x)=-x^{3/2}$.
\end{enumerate}
In both these examples  there is insufficient smoothness in the function $q$ for Lohner's code to compute 
enclosures at $x=0$. We therefore compute enclosures at $x=\epsilon$ and use the above theorem 
to obtain enclosures at $x=0$.
\par
In the first example we choose $f$ and $g$ to be the first-order Taylor polynomial approximations to $y$ and  $y^{'}$ expanded about $x=\epsilon$. Writing $A_0=y(\epsilon)$ and  $ A_1=y^{'}(\epsilon)$ these are
respectively
\begin{eqnarray*}
f (x)&=& A_0 + ( x-\epsilon)A_1, \\
g(x) &=& A_1 + (x-\epsilon) (- \sqrt{\epsilon}-\lambda)A_0.
\end{eqnarray*} 
Since for this example
\begin{eqnarray*}
\int_0^\epsilon \mid c \mid dt  &\leq & \frac{2}{3}\epsilon ^{3/2} + \mid \lambda \mid \epsilon \\
\int_0^\epsilon \mid c -c(\epsilon)\mid dt  &= & \frac{1}{3} \epsilon ^{3/2} \\
\int_0^\epsilon \mid c \mid ( \epsilon -t)dt  &\leq & \frac{4}{15}\epsilon ^{5/2} +\frac{1}{2} \mid \lambda \mid \epsilon^2,
\end{eqnarray*}
the inequalities (\ref{eq:5.1}) and  (\ref{eq:5.2}) become
\begin{eqnarray*}
 \mid y(0)-(A_0 - \epsilon A_1) &\mid&   \leq 
\frac{\epsilon^2}{1-\epsilon^2 ( \mid \lambda \mid + \frac{2}{3} \sqrt{\epsilon} ) } \times
\nonumber \\
&& \left \{  (  \frac{1}{2} \mid \lambda \mid  + \frac{5}{6} \sqrt{\epsilon} ) \mid A_0 \mid + \epsilon ( 
\frac{1}{2} \mid \lambda \mid +\frac{4}{15}\sqrt{\epsilon }) \mid A_1 \mid 
\right \} \nonumber \\
\mid y^{'}(0)-(A_1 + \epsilon( \lambda + \sqrt{\epsilon} ) A_0) &\mid&   \leq   
\frac{\epsilon^{3/2}}{1-\epsilon^2 ( \mid \lambda \mid +\frac{2}{3} \sqrt{\epsilon} ) } \times
\nonumber \\
&&
\left \{
  [  \frac{1}{2} \epsilon^{3/2}( \mid \lambda \mid  + \sqrt{\epsilon} ) ( \mid \lambda \mid + \frac{2}{3}\sqrt{\epsilon}) + \frac{1}{3
} ] \mid A_0 \mid + \sqrt{\epsilon} ( 
\frac{1}{2} \mid \lambda \mid +\frac{4}{15} \sqrt{\epsilon }) \mid A_1 \mid 
\right \} \label{eq:5.3}
\end{eqnarray*} 
In the example that we report on, $\lambda=1+i$ and we have taken $\epsilon = 0.000015625$. However,  since   Lohner's
code uses a fixed step size algorithm, in order to 
 achieve such a small value of $\epsilon$ we have integrated over $[10,0.03125]$ with a step size of $0.03125$, then reduced the step size to $0.000015625$ to perform the integration over $[0.03125,0.000015625]$.
This yields an enclosure 
\begin{displaymath}
m(1+i) = 0.25936_{61}^{82} + 0.6719_{78}^{82} i.
\end{displaymath}
We have compared our result with  the estimate of  $m(1+i)$ obtained from the Brown Kirby Pryce code \cite{BrownKirbyPryce1}. That algorithm gives $m(1+i)=0.25937860 + 0.67196464 i$ which is slightly outside our enclosures and reflects the lack of smoothness in  $q$ experienced by the Runga-Kutta method they employ. 
\par
In the next example we take $q=-x^{3/2}$ and as before  we take $f$ and $g$ to be the first-order 
Taylor approximations  of $y$ and $y^{'}$ respectively, expanded about $x=\epsilon$.
This time  we get
 \begin{eqnarray*}
 \mid y(0)-(A_0 - \epsilon A_1) \mid   &\leq& \frac{\epsilon^2}{1- \epsilon^{2}(\mid \lambda \mid+\frac{2}{5}\epsilon^{3/2})  }   \\
&& \{
(\frac{1}{2} \mid \lambda \mid +\frac{11}{10} \epsilon^{3/2} ) \mid A_0 \mid + \epsilon (\frac{1}{2} \mid \lambda \mid + \frac{4}{35}
\epsilon^{3/2}) \mid A_1 \mid
  \}
 \\
\mid y^{'}(0)-(A_1 + \epsilon( \lambda + \epsilon^{3/2} ) A_0) \mid   &\leq&\frac{\epsilon^2}{1-\epsilon^{2}(\mid \lambda \mid +\frac{2}{5}\epsilon^{3/2}) } \\
&&\{
 [\frac{1}{2} \epsilon (  \mid  \lambda \mid  + \epsilon ^{3/2} )(     \mid \lambda \mid  + \frac{2}{5} \epsilon^{3/2}) +\frac{3}{5}\sqrt{ \epsilon}]\mid A_0 \mid 
+  (\frac{1}{2} \mid \lambda \mid + \frac{4}{35} \epsilon^{3/2})   \mid A_1 \mid \}.
\end{eqnarray*}
This gives an enclosure
\begin{displaymath}
m(1+i) = 0.345229_{19}^{40} + 0.705227_{69}^{85} i.
\end{displaymath}
\par
We have investigated the possibility of choosing $f$ and $g$ to be the second-order Taylor polynomials. However, for both the examples that we  have considered there appears to be no   improvement in the bound over that achieved by linear functions $f,g$.

\scriptsize
\bibliographystyle{unsrt}
\bibliography{../bib_dir/bib1,../bib_dir/bibliography,../bib_dir/help,../bib_dir/manual}
\normalsize

\section*{Author Addresses}
   (BMB), (MSPE), (DKRM) - Department of Computer Science, Cardiff University of Wales, Cardiff, CF2 3XF, U.K, 
\newline
(MP) - Mathematisches Institut I, Universit\"at Karlsruhe,76128  Karlsruhe, Germany

\section*{Date}

Eingegamgen am 26. November 1996

\end{document}